\documentclass[11pt,reqno]{article}
\usepackage{jmlr2e}
\usepackage{fullpage,times,graphicx,amssymb,amsmath,amsfonts,bbm,psfrag,xcolor}

\usepackage{graphicx}
\usepackage{amssymb}
\usepackage{amsmath}
\usepackage{amsfonts}
\usepackage{bbm}
\usepackage{psfrag}
\usepackage{xcolor}
\usepackage{fullpage}
\usepackage{epstopdf}
\usepackage{graphicx}
\usepackage{times}
\usepackage{jmlr2e}
\usepackage{commath}
\usepackage{thmtools}
\usepackage{thm-restate}
\usepackage{times}
\usepackage{graphicx}
\usepackage{natbib}

\definecolor{ddarkbrown}{rgb}{0.5,0.2,0.05} \definecolor{bbluegray}{rgb}{0.05,0,0.5}

\usepackage[colorlinks,citecolor=bbluegray,linkcolor=ddarkbrown,urlcolor=blue,breaklinks]{hyperref}

\usepackage{dirtytalk}
\usepackage{subfig,float} 

\usepackage{mathtools}

\usepackage{algorithm,algcompatible,algpseudocode}
\usepackage{multicol,lipsum}
\algnewcommand{\Inputs}[1]{%
	\State \textbf{Inputs: \:}{#1}
}

\algnewcommand{\Output}[1]{%
	\State \textbf{Output: \:}{#1}
}
\algnewcommand{\Initialize}[1]{%
	\State \textbf{Initialize: \:}{#1}
}

\algnewcommand{\IIf}[1]{\State\algorithmicif\ #1\ \algorithmicthen}
\algnewcommand{\EndIIf}{\unskip\ \algorithmicend\ \algorithmicif}

\let \oldsection \section
\renewcommand{\section}{\vspace{3ex plus 1ex}\oldsection}

\newcommand{\BEAS}{\begin{eqnarray*}}
	\newcommand{\EEAS}{\end{eqnarray*}}
\newcommand{\BEA}{\begin{eqnarray}}
\newcommand{\EEA}{\end{eqnarray}}

\newcommand{\BEQ}{\begin{equation}}
\newcommand{\EEQ}{\end{equation}}
\newcommand{\BIT}{\begin{itemize}}
	\newcommand{\EIT}{\end{itemize}}
\newcommand{\BNUM}{\begin{enumerate}}
	\newcommand{\ENUM}{\end{enumerate}}

	\newcommand{\D}{\mathcal{D}}

				\newcommand{\G}{\mathcal{G}}

	\newcommand{\R}{\mathbb{R}}

 \newcommand{\M}{\mathcal{M}}

\newcommand{\Proj}{\mathcal{P}}
\newcommand{\T}{\mathcal{T}}

\newcommand{\K}{\mathcal{K}}

\newcommand{\BA}{\begin{array}}
	\newcommand{\EA}{\end{array}}

 \numberwithin{dummy}{section}

\numberwithin{mythm}{section}
\numberwithin{mydef}{section}
\numberwithin{myprop}{section}
\numberwithin{mylem}{section}
\numberwithin{mycor}{section}

\title{Equivariance Regularization for Image Reconstruction}

\begin{document}
	\author{Junqi Tang \email jt814@cam.ac.uk\\
		\addr Department of Applied Mathematics and Theoretical Physics (DAMTP),\\ University of Cambridge
	}
	\editor{}
	
	\maketitle

%\bibliography{a_note.bib}
%\bibliographystyle{agsm}
%\appendix

\begin{abstract}
In this work, we propose Regularization-by-Equivariance (REV), a novel structure-adaptive regularization scheme for solving imaging inverse problems under incomplete measurements. This regularization scheme utilizes the equivariant structure in the physics of the measurements -- which is prevalent in many inverse problems such as tomographic image reconstruction -- to mitigate the ill-poseness of the inverse problem. Our proposed scheme can be applied in a plug-and-play manner alongside with any classic first-order optimization algorithm such as the accelerated gradient descent/FISTA for simplicity and fast convergence. The numerical experiments in sparse-view X-ray CT image reconstruction tasks demonstrate the effectiveness of our approach.
\end{abstract}

\section{Introduction}

Imaging from incomplete measurements has become an important yet challenging problem studied intensely by researchers throughout recent decades. Such imaging systems, for example the sparse-view X-ray CT and compressed-sensing MRI, can be generally expressed as:
\begin{equation}
    b = Ax^\dagger + w,
\end{equation}
where $x^\dagger \in \mathbb{R}^d$ denotes the (vectorized) ground-truth image to be inferred, $A \in \mathbb{R}^{n \times d}$ denotes the forward measurement operator, $w \in \mathbb{R}^n$ denotes the noise (could be measurement-dependent), and $b \in \mathbb{R}^d$ denotes the measurement data. 

Here we are mostly interested in the case where $\mathrm{rank}(A) < d$. To obtain a high quality solution, one must cope with the non-trivial null space of the measurement operator. The traditional approach is to solve a composite convex program of the form:
\begin{equation}\label{obj1}
    x^\star \in \arg\min_{x} f(Ax, b) + \lambda r(x),
\end{equation}
where the first term is the data-consistency loss, with one typical choice being the least-squares fit  $f(Ax, b) = \frac{1}{2n} \|Ax - b\|_2^2$, while the second term is a convex regularization modeling the image prior. Classical choices of $r(x)$ includes $\ell_1$ regularization in the wavelet domain, total-variation (TV) regularization and its variants \citep{chambolle2016introduction}, etc. Such choices of $r(x)$ encode certain structures of the images, for example sparsity in some transformed domains. In the literature of statistical inference, researchers \citep{agarwal2009information,agarwal2010fast,agarwal2012fast, 2015_Pilanci_Randomized,oymak2017sharp,pmlr-v70-tang17a,qu2016linear, tang2018rest,tang2018structure} have shown theoretical robust recovery results for (\ref{obj1}) under the case where $A$ is a random matrix drawn from a Gaussian-type distribution. However, most real-world inverse problems have a deterministic measurement forward operator which is highly structured, and in such cases there is no non-trivial recovery guarantee. We believe that these classical regularization schemes, as well as recently occurred plug-and-play/regularization-by-denoising \citep{venkatakrishnan2013plug,romano2017little, reehorst2018regularization} with advanced denoisers \citep{dabov2007image,buades2005non} and deep-learning based priors \citep{zhang2017beyond, chen2017trainable, tachella2020neural, jin2017deep,ulyanov2018deep}, have a fundamental limitation, that they do not have any explicit mechanism to take into account the structure of the forward measurement operator $A$.

Inspired by recent works in computer vision and deep learning literature \citep{celledoni2021equivariant,chen2021equivariant}, we propose a framework for instrumenting structure-adaptive regularization, which utilizes the prevalent equivariant nature of the measurement systems in imaging applications.

\section{Regularization-by-Equivariance (REV)}
In this section we present our Regularization-by-Equivariance (REV) framework, which is tailored to the inherent equivariant structure of the inverse problems.
\subsection{Algorithmic Framework of REV}
We propose to solve the following optimization program with our Regularization-by-Equivariance (REV) scheme:
\begin{equation}\label{rev_obj}
     x^\star \in \arg\min_{x \in \K} f(Ax, b) + \lambda \mathbb{E}_{g} [x^T(x - T_g^* \D( T_g x))],
\end{equation}
where we denote $\D$ as a user-defined denoising/artifact removal algorithm, for example the BM3D \citep{dabov2007image} or a pre-trained\footnote{If we use a untrained network we would need to run a DIP-type of iterations \citep{ulyanov2018deep} with stochastic alternating minimization \citep{driggs2020spring}.} neural-network mapping \citep{zhang2017beyond}, etc. Alternatively we may also choose $\D$ to explicitly encode the forward map $\D(\T_g x) = \M(A \T_g x))$, where $\M$ is a pretrained inverse map of $A$. We denote here $g \in \G$ where $\G$ is a group of transformations, and $T_g \in \R^{d \times d}$ are unitary matrices such that for a certain subset $X \subset \R^d$ containing real-world/clinical tomographic images:
\begin{equation}
    T_g x \in X,\ \forall x \in X.
\end{equation}
For sparse-view CT, a typical choice of $\G$ would be random rotations. Meanwhile $\K$ is a constraint set, for example we always use at least the box-constraint to restrict the pixel values within a certain range in practice. %If further the $\M$ we choose is a linear operator such as FBP, the optimization problem (\ref{rev_obj}) is convex.

Our REV is related to the regularization-by-denoising (RED) \citep{romano2017little,reehorst2018regularization}, where they select beforehand a denoiser $\D(\cdot)$, for example the BM3D \citep{dabov2007image} or the DnCNN \citep{zhang2017beyond}, and set the regularizer as $r(x) = x^T(x - \D(x))$. A practical choice of approximate the gradient of RED is simply $x - \D(x)$ (if $\D$ satisfies certain conditions, this approximation is exact). Inspired by RED, we also propose to approximate the gradient of REV term $\mathbb{E}_{g} [x^T(x - T_g^* \D( T_g x))]$ by first randomly sample a operation $T_g$ from $\G$, and approximate the gradient as:
\begin{equation}
    \nabla_x \mathbb{E}_{g} [x^T(x - T_g^* \D( T_g x))] \approx x - T_g^* \D (T_g x).
\end{equation}
Now we present our algorithmic framework\footnote{In practice very often the transformation $T_g$ can only be performed approximately (we denote this approximation as $\T_g$) using interpolation for FISTA-REV, due to the fact that we are working on discretized image domains.} with accelerated-gradient/FISTA \citep{beck2009fast,nesterov2007gradient,chambolle2015convergence} iterations (we could also choose stochastic gradients \citep{johnson2013accelerating,allen2017katyusha,tang2018rest} instead on the data-fit\footnote{As discussed in \citep{tang2020practicality,tang2019limitation}, in some inverse problems such as X-ray CT, the stochastic gradient methods can provide much improved convergence rates over deterministic gradient methods. However, in some other scenarios such as compressive sensing MRI and space-variant deblurring, the stochastic gradient methods are not efficient.}):

 \begin{eqnarray*}
 && \mathrm{\textbf{FISTA-REV}} - \mathrm{Initialize}\ x_{0}\in \R^d, \ y_0 \in \R^d, \ a_0 = 1\\
 &&\mathrm{For} \ \ \ k = 0, 1, 2,...,  K\\
&&\left\lfloor
\begin{array}{l}
\mathrm{Sample}\ g \in \G\\
x_{k+1} =  \Proj[y_k - \eta \cdot [ \nabla_{y_k} f(Ay_k, b) + \lambda(y_k - \T_g^* \D( \T_g y_k))]]\\
a_{k+1} = (1 + \sqrt{1 + 4a_k^2})/2;\\
y_{k + 1} = x_{k+1} + \frac{a_{k} - 1}{a_{k+1}} (x_{k+1} - x_{k}).
\end{array}
\right.
 \end{eqnarray*}
 
\footnote{We may also sample a number of $g$ in each iteration and compute the average of REV gradients to reduce the estimation variance \citep{defazio2014saga}, and we numerically find this would further improve the convergence rates and reconstruction quality.}In practice, very often we cannot exactly perform the group actions in a finite-dimensional space, for example, the rotation operations with arbitrary angles have to be performed approximately, via interpolation tricks. Here we redefine such an approximated operation using interpolation as $\T_g$.

Here we often choose $\Proj$ as the projection to some constraint set $\K$ (for example a box constraint for restricting pixel-values). Actually the choices of $\Proj$ and $\D$ can be very flexible, such as the proximal operator of some extra regularizer, or a denoising algorithm (BM3D/NLM/TNRD \citep{chen2017trainable}), or a pretrained neural-denoiser (DnCNN). We can also consider a deep-unrolling \citep{adler2018learned, tang2021stochastic} of FISTA-REV, where we would unfold the FISTA-REV iterations and learn both $\Proj_k$ and $\M_k$ for each iteration from training-data. Such a deep unrolling network will have a built-in element of equivariance.

\subsection{A simplified scheme utilizing the continuous image domain}

In this subsection, we also propose a simple reduced scheme of FISTA-REV which we will show in the experiment section having excellent numerical performance (despite being a much simplified scheme).

If just for now we choose to explicitly encode the forward map and assume that we have a very high-quality inverse map $\M$ of $A$ somehow, the following approximation holds:
\begin{equation}
    \nabla_x \mathbb{E}_{g} [x^T(x - \T_g^* \D( \T_g x))] \approx x - \T_g^* \M(A \T_g x) \approx x -  \T_g^*  \T_g x,
\end{equation}
where $ \T_g^*  \T_g \neq I$, due to the fact that $\T_g$ is an approximated action with interpolation, utilizing the continuous image domain. This idealized thinking process leads us to a very simple instance of REV:

 \begin{eqnarray*}
 && \mathrm{\textbf{Simplified FISTA-REV}} - \mathrm{Initialize}\ x_{0}\in \R^d, \ y_0 \in \R^d, \ a_0 = 1\\
 &&\mathrm{For} \ \ \ k = 0, 1, 2,...,  K\\
&&\left\lfloor
\begin{array}{l}
\mathrm{Sample}\ g \in \G\\
x_{k+1} =  \Proj[y_k - \eta \cdot [ \nabla_{y_k} f(Ay_k, b) + \lambda(y_k - \T_g^*  \T_g y_k)]]\\
a_{k+1} = (1 + \sqrt{1 + 4a_k^2})/2;\\
y_{k + 1} = x_{k+1} + \frac{a_{k} - 1}{a_{k+1}} (x_{k+1} - x_{k}).
\end{array}
\right.
 \end{eqnarray*}
 
Although this simplified version of REV does not explicitly involve the measurement operator and denoiser, we find that if the underlying interpolation we choose is appropriately aligned with the null-space of $A$, this simplified REV will still provide powerful regularization to the inverse problem. For the sparse-view CT, if we choose $\Proj$ to be a projection towards a box-constraint, and $\T_g$ to be the rotations with linear interpolations, the Simplified FISTA-REV solves a convex optimization problem.

 \begin{figure}[t]
   \centering

    {\includegraphics[width= .975\textwidth]{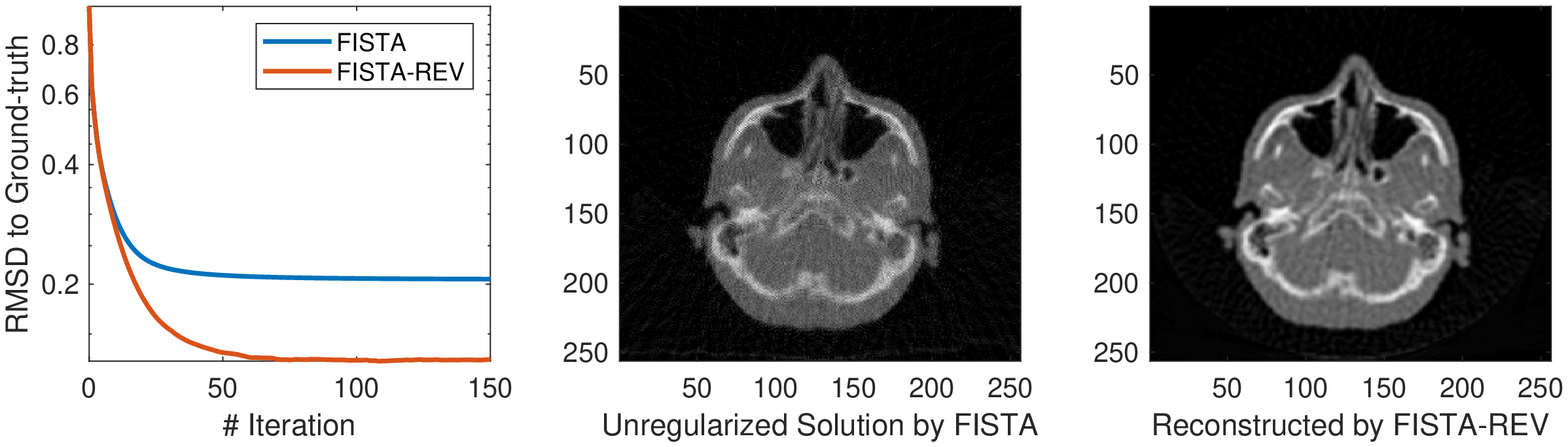}}
    {\includegraphics[width= .975\textwidth]{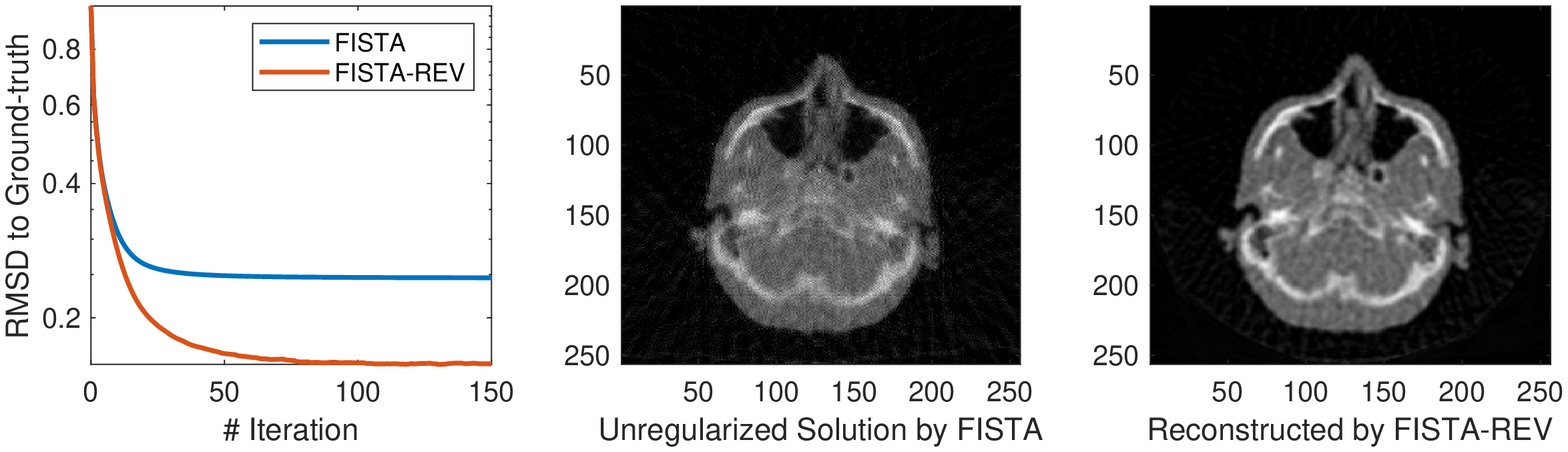}}
    {\includegraphics[width= .975\textwidth]{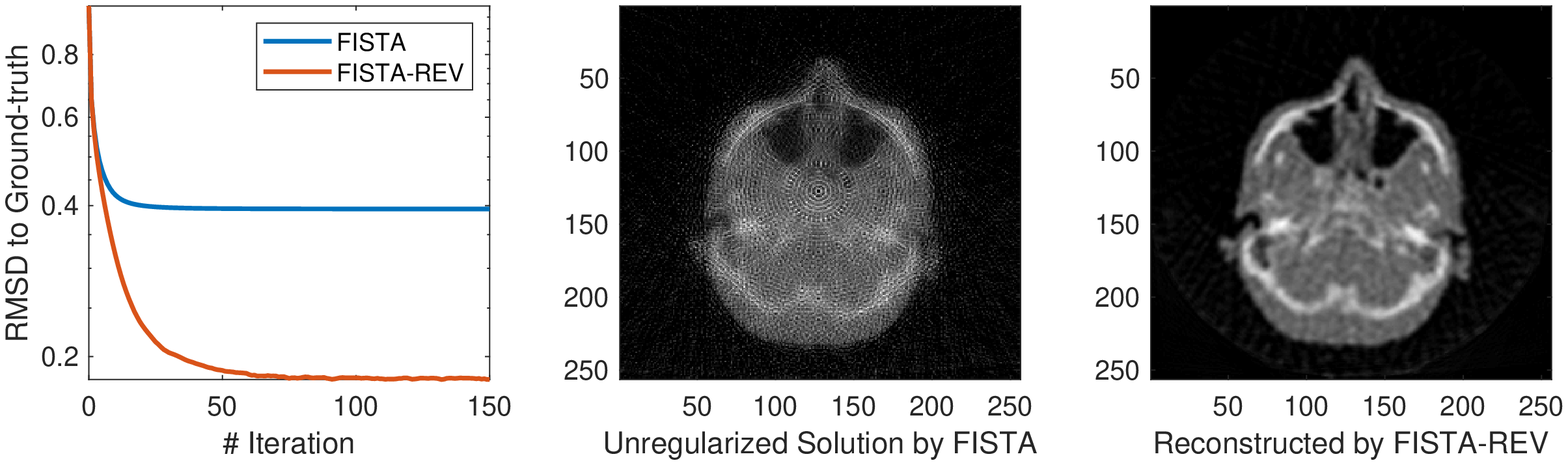}}

   \caption{Sparse-view fan-beam CT (with $I_0 = 2 \times 10^{7.5}$) reconstruction results for FISTA and FISTA-REV, without extra regularization. The first sparse view CT example (first-row): $A_1 \in \R^{13650 \times 65536}$, where the number of views is 60. The second example: $A_2 \in \R^{9120 \times 65536}$, where the number of views is 40. The third sparse view CT example we consider has a forward operator $A_3 \in \R^{4560 \times 65536}$, where the number of views is 40. The third example is extremely sparse-view.}
    \label{rev_e1}
\end{figure}

 \begin{figure}[t]
   \centering

    {\includegraphics[width= .975\textwidth]{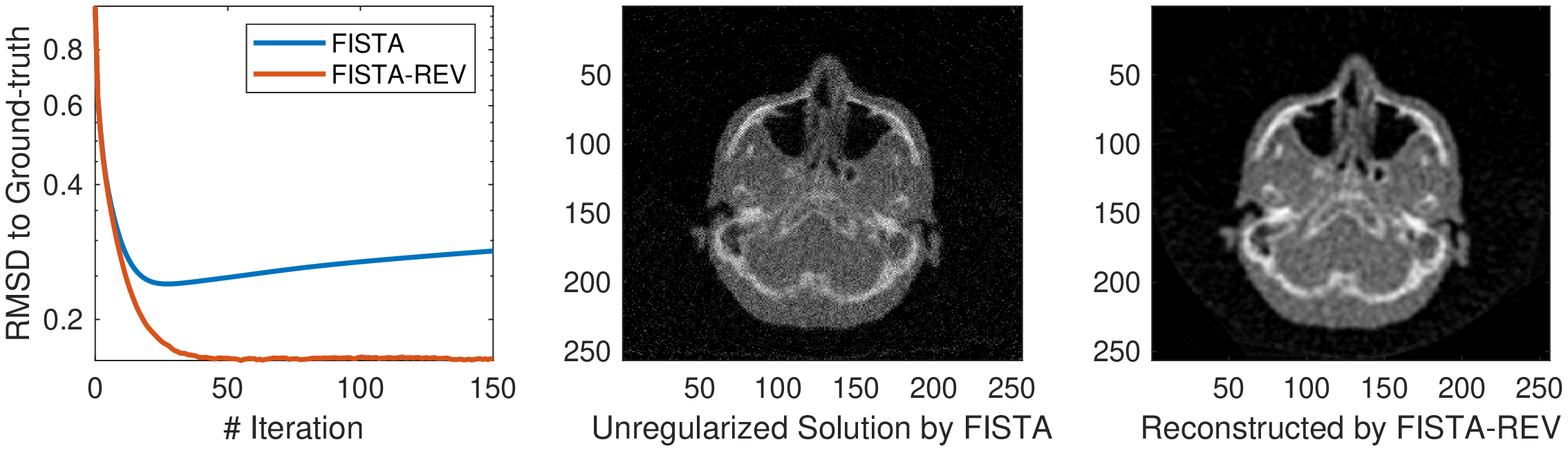}}
    {\includegraphics[width= .975\textwidth]{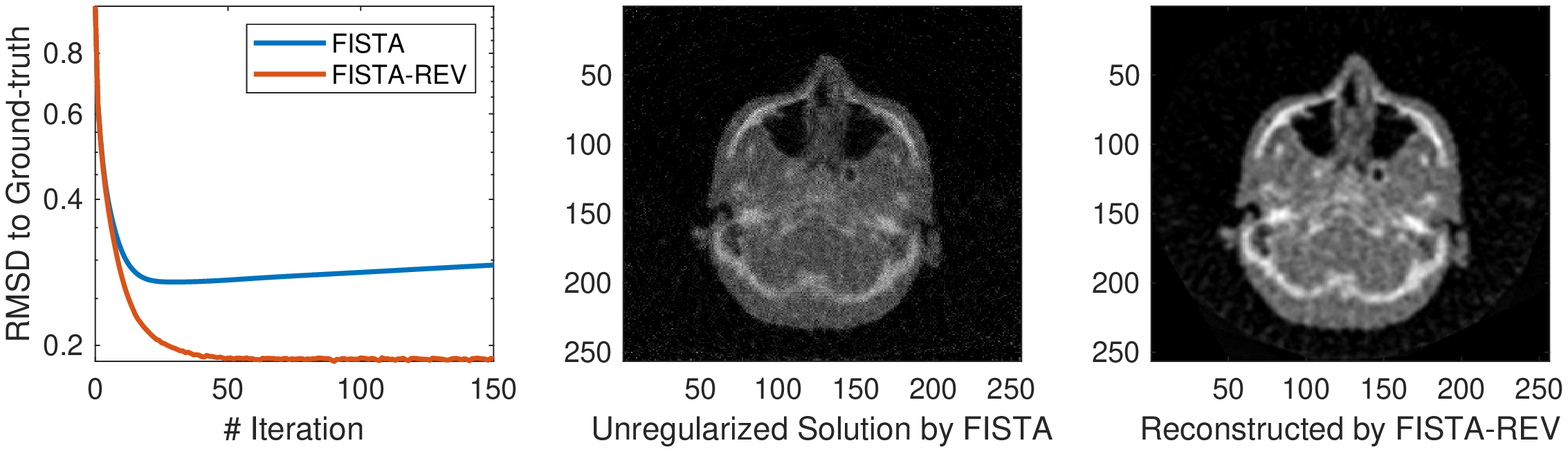}}
    {\includegraphics[width= .975\textwidth]{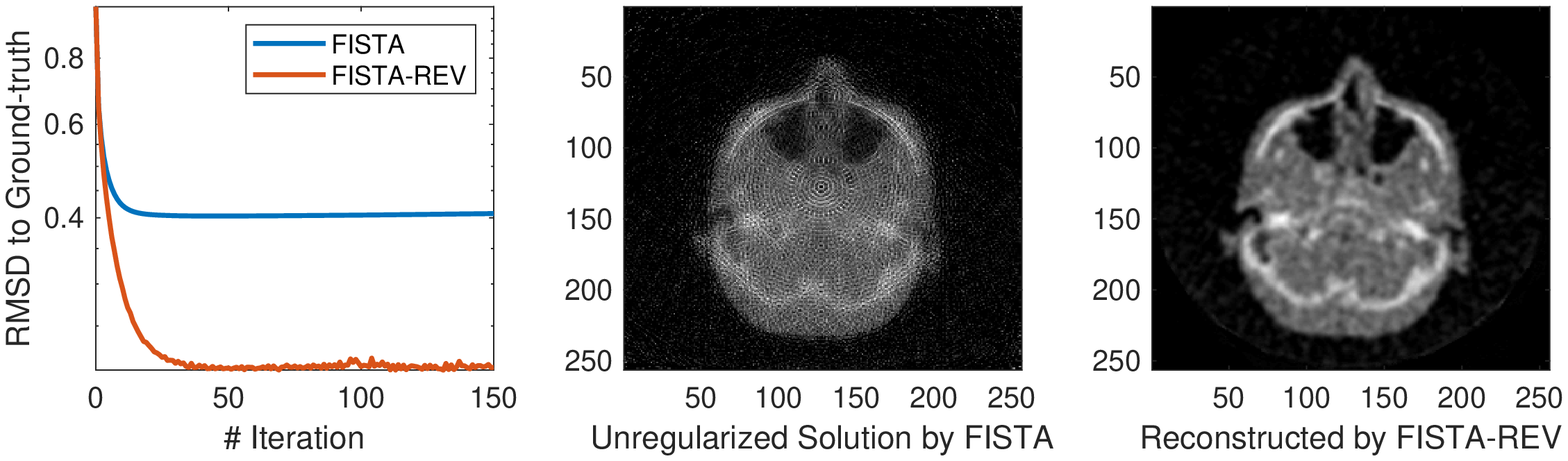}}

   \caption{Sparse-view low-dose fan-beam CT (with $I_0 = 2 \times 10^{3.5}$) reconstruction results for FISTA and FISTA-REV, without extra regularization. The first sparse view CT example (first-row): $A_1 \in \R^{13650 \times 65536}$, where the number of views is 60. The second example: $A_2 \in \R^{9120 \times 65536}$, where the number of views is 40. The third sparse view CT example we consider has a forward operator $A_3 \in \R^{4560 \times 65536}$, where the number of views is 40. The third example is extremely sparse-view.}
    \label{rev_e2}
\end{figure}

\section{Numerical Experiments}

In this section we perform numerical experiments on sparse-view CT. We start by the Simplied FISTA-REV without any extra image priors, for a proof-of-concept example. Then we will present the comparison between REV and RED \citep{romano2017little,reehorst2018regularization}, where we choose BM3D as the denoiser for both REV and RED.

\subsection{Motivational Experiments on the Simplified FISTA-REV}

For the proof of concept, we start by presenting here some preliminary results of (simplified) FISTA-REV on sparse-view X-ray CT, where we set the measurements to be massively incomplete. Here we run all the experiments in MATLAB R2018a. We choose $\T_g$ to be rotations with random angles $\theta \in (0, 360)$, while $\T_g^*$ to be rotation with angle $-\theta$, both implemented using MATLAB's \texttt{imrotate} function, and we choose the bicubic interpolation method.

Utilizing the Beer-Lambert law, we simulate the projection data of the sparse-view fan-beam CT, corrupted with Poisson noise:
\begin{equation}
    b \sim \mathrm{Poisson}(I_0 e^{-Ax^\dagger}),
\end{equation}
and we take the logarithmic to linearize the measuements.

 \begin{figure}[t]
   \centering

    {\includegraphics[width= .975\textwidth]{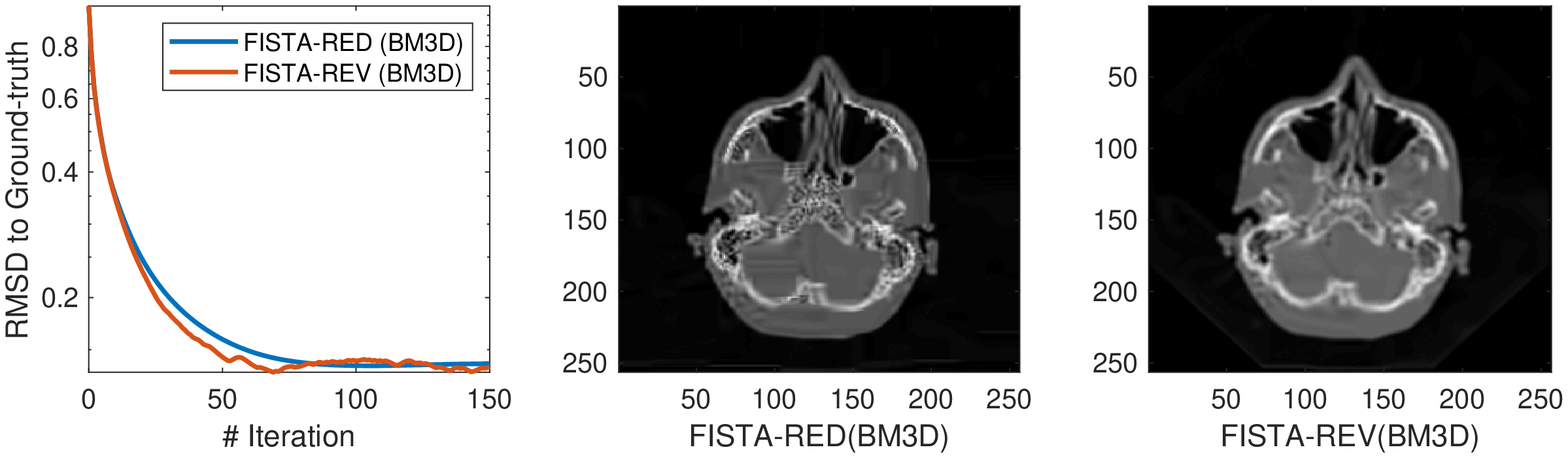}}
    {\includegraphics[width= .975\textwidth]{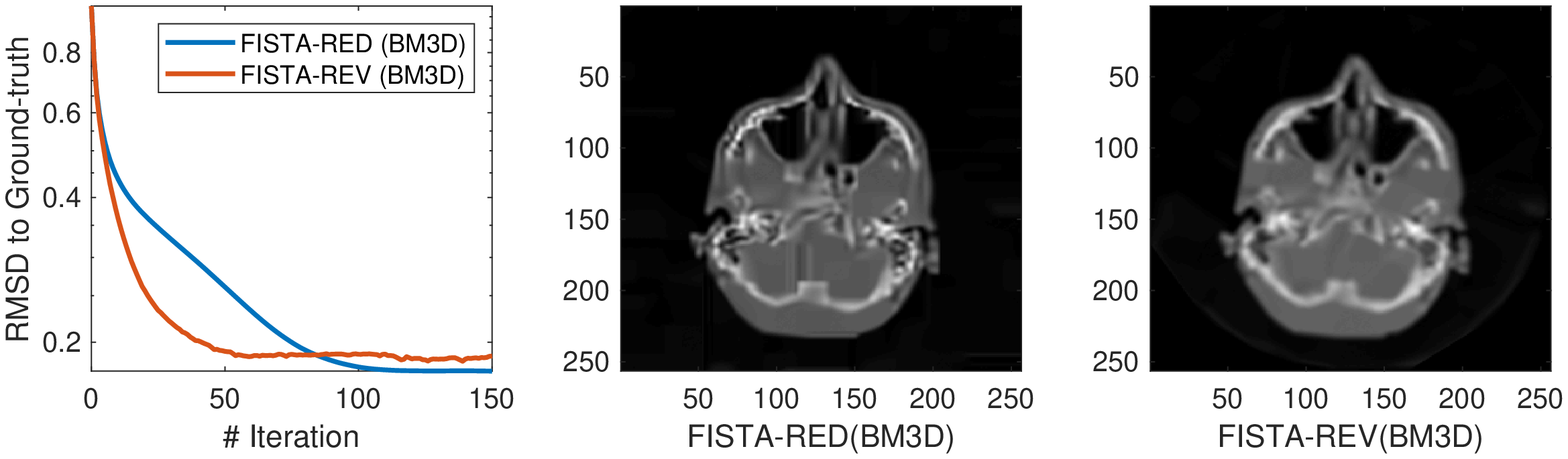}}
    {\includegraphics[width= .975\textwidth]{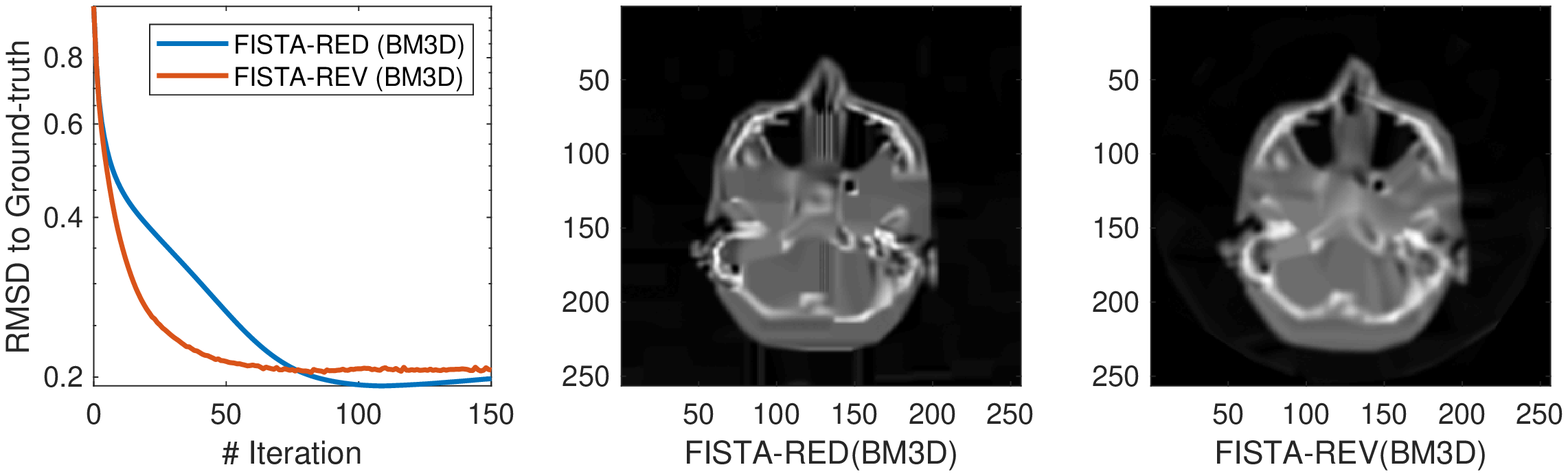}}

   \caption{Sparse-view low-dose fan-beam CT (with $I_0 = 2 \times 10^{3.5}$) reconstruction results for FISTA and FISTA-REV {\color{purple}(with BM3D denoisers)}. The first sparse view CT example (first-row): $A_1 \in \R^{13650 \times 65536}$, where the number of views is 60. The second example: $A_2 \in \R^{4560 \times 65536}$, where the number of views is 40. The third sparse view CT example we consider has a forward operator $A_3 \in \R^{3420 \times 65536}$, where the number of views is 30. The second and third example is extremely sparse-view.}
    \label{rev_e3}
\end{figure}

We consider 3 highly-underdetermined sparse-view CT problems. The ground-truth $x^\dagger \in \R^{65536}$ (sized $256 \times 256$). The first sparse view CT example we consider has a forward operator $A_1 \in \R^{13650 \times 65536}$, where the number of views is 60 (equally spaced), each view has 228 measurements. The second sparse view CT example we consider has a forward operator $A_2 \in \R^{9120 \times 65536}$, where the number of views is 40, each view has 228 measurements. The third sparse view CT example we consider has a forward operator $A_3 \in \R^{4560 \times 65536}$, where the number of views is 40, each view has only 114 measurements, making it extremely under-determined. We set $\Proj$ to be the projection operator on a box constraint restricting the pixel values to live within the range of $[0, 1]$ for both FISTA and FISTA-REV. For clarity of the comparison here we do not use any extra regularization.

In Figure \ref{rev_e1} and \ref{rev_e2}, we report the root mean square distance (RMSD, which is computed as: $\frac{\|x - x^\dagger\|_2}{\sqrt{d}}$) of the iterates of FISTA and FISTA-REV towards the ground-truth, and also the outcome images of the algorithms at termination. The first row presents the results for the first example, and so on. We can observe that in these massively under-determined examples, FISTA-REV consistently provides an improvement over FISTA by a large margin, demonstrating the effectiveness of our regularization scheme which is tailored to mitigate the null-space artifacts. 

Meanwhile, noting that the number of measurements decreases from task 1 to task 3. If we compare the RMSD plots and images vertically, we can observe that the performance of FISTA-REV degrades much slower than the unregularized FISTA solution as the number of measurements decreases.

Most remarkably, in the third example which is an extremely sparse-view case, we can observe that the unregularized solution of FISTA has a massive amount of artifacts due to the huge null-space, while our FISTA-REV still provides reasonable reconstruction performance. From this extreme example we can clearly see that the null-space aritifacts occur as circle-like stripping shapes, which can be indeed massively mitigated using random rotation group actions with interpolation tricks.

\subsection{REV versus RED}

Now we turn to the comparison of REV with the Regularization-by-Denoising (RED, \cite{romano2017little}) which corresponds to REV without any \say{null-space-compensating transforms} ($\T_g$). We use the BM3D \citep{dabov2007image} denoiser as $\D$ for both REV and RED.

Here we also consider 3 highly-underdetermined sparse-view CT problems. The first sparse view CT example we consider has the same forward operator $A_1 \in \R^{13650 \times 65536}$ as the first example in the previous section. The second sparse view CT example we consider has a forward operator $A_2 \in \R^{4560 \times 65536}$, where the number of views is 40, each view has only 114 measurements. The third sparse view CT example we consider has a forward operator $A_3 \in \R^{3420 \times 65536}$, where the number of views is 30, each view has also only 114 measurements. The second and third examples are extremely under-determined. We simulate low-dose measurements with $I_0 = 2 \times 10^{3.5}$ for reconstruction.

We present the results of FISTA-RED and FISTA-REV in Figure \ref{rev_e3}. We can observe that, although FISTA-REV and FISTA-RED converges towards solutions with similar RMSD, the FISTA-RED solutions has obvious artifacts (likely due to its non-adaptivity to the large-scale null-space), while FISTA-REV is free of these artifacts. 

Interestingly, we can also observe that for the extreme sparse-view cases, the FISTA-REV converges much faster than FISTA-RED \footnote{We choose the same step-sizes and regularization parameters for both of the algorithms. The only difference between the two algorithms is with/without rotations.}, suggesting that the underlying inverse problem became effectively better-conditioned due to the rotation transformations. This effect again clearly demonstrates the importance and effectiveness of constructing null-space structure-adaptive regularization.

\small
\bibliography{main.bib}

\end{document}